\documentclass[12pt]{article}

\usepackage{amsmath,amssymb}

\begin{document}

\title{A semi-local trace identity and the Riemann hypothesis for function fields}
\author{Anton Deitmar}

\date{}
\maketitle

{\bf Abstract} The asymptotic trace formula of Connes is
restated in a semi-local form, thus showing that the
difficulties in proving it directly do not lie in the
change of topology when transgressing from finitely many
to infinitely many places.

\tableofcontents

\def \1{{\bf 1}}
\def \a{{{\mathfrak a}}}
\def \ad{{\rm ad}}
\def \al{\alpha}
\def \ar{{\alpha_r}}
\def \A{{\mathbb A}}
\def \Ad{{\rm Ad}}
\def \Aut{{\rm Aut}}
\def \b{{{\mathfrak b}}}
\def \bs{\backslash}
\def \B{{\cal B}}
\def \c{{\mathfrak c}}
\def \cent{{\rm cent}}
\def \C{{\mathbb C}}
\def \CA{{\cal A}}
\def \CB{{\cal B}}
\def \CC{{\cal C}}
\def \CE{{\cal E}}
\def \CF{{\cal F}}
\def \CG{{\cal G}}
\def \CH{{\cal H}}
\def \CHC{{\cal HC}}
\def \CL{{\cal L}}
\def \CM{{\cal M}}
\def \CN{{\cal N}}
\def \CP{{\cal P}}
\def \CQ{{\cal Q}}
\def \CO{{\cal O}}
\def \CR{{\cal R}}
\def \CS{{\cal S}}
\def \CT{{\cal T}}
\def \CV{{\cal V}}
\def \det{{\rm det}}
\def \diag{{\rm diag}}
\def \dist{{\rm dist}}
\def \eqn{\begin{eqnarray*}}
\def \endeqn{\end{eqnarray*}}
\def \End{{\rm End}}
\def \F{{\mathbb F}}
\def \Fx{{\mathfrak x}}
\def \FX{{\mathfrak X}}
\def \g{{{\mathfrak g}}}
\def \ga{\gamma}
\def \Ga{\Gamma}
\def \h{{{\mathfrak h}}}
\def \Hom{{\rm Hom}}
\def \im{{\rm im}}
\def \Im{{\rm Im}}
\def \Ind{{\rm Ind}}
\def \k{{{\mathfrak k}}}
\def \K{{\cal K}}
\def \l{{\mathfrak l}}
\def \la{\lambda}
\def \lap{\triangle}
\def \La{\Lambda}
\def \m{{{\mathfrak m}}}
\def \mod{{\rm mod}}
\def \n{{{\mathfrak n}}}
\def \name{\bf}
\def \N{\mathbb N}
\def \o{{\mathfrak o}}
\def \ord{{\rm ord}}
\def \O{{\cal O}}
\def \p{{{\mathfrak p}}}
\def \ph{\varphi}
\def \prf{\noindent{\bf Proof: }}
\def \Per{{\rm Per}}
\def \q{{\mathfrak q}}
\def \qed{\ifmmode\eqno Q.E.D.\else\noproof\vskip 12pt plus 3pt minus 9pt \fi}
 \def\noproof{{\unskip\nobreak\hfill\penalty50\hskip2em\hbox{}%
     \nobreak\hfill Q.E.D.\parfillskip=0pt%
     \finalhyphendemerits=0\par}}
\def \Q{\mathbb Q}
\def \res{{\rm res}}
\def \R{{\mathbb R}}
\def \Re{{\rm Re \hspace{1pt}}}
\def \r{{\mathfrak r}}
\def \ra{\rightarrow}
\def \rank{{\rm rank}}
\def \supp{{\rm supp}}
\def \Spin{{\rm Spin}}
\def \t{{{\mathfrak t}}}
\def \T{{\mathbb T}}
\def \tr{{\hspace{1pt}\rm tr\hspace{2pt}}}
\def \vol{{\rm vol}}
\def \z{\zeta}
\def \Z{\mathbb Z}
\def \={\ =\ }

\newcommand{\frack}[2]{\genfrac{}{}{0pt}{}{#1}{#2}}
\newcommand{\rez}[1]{\frac{1}{#1}}
\newcommand{\der}[1]{\frac{\partial}{\partial #1}}
\newcommand{\norm}[1]{\parallel #1 \parallel}
\renewcommand{\matrix}[4]{\left( \begin{array}{cc}#1 & #2 \\ #3 & #4 \end{array}
            \right)}
\renewcommand{\sp}[2]{\langle #1\,, #2\rangle}

\newtheorem{theorem}{Theorem}[section]
\newtheorem{conjecture}[theorem]{Conjecture}
\newtheorem{lemma}[theorem]{Lemma}
\newtheorem{corollary}[theorem]{Corollary}
\newtheorem{proposition}[theorem]{Proposition}

$$ $$

\begin{center} {\bf Introduction} \end{center}

The Riemann hypothesis for Hecke $L$-series ``mit Gr\"o\ss
encharakteren'' over function fields has been proved by
Weil \cite{weil3}. In \cite{bomb}, Bombieri gave a
simplified proof that is based on the Riemann-Roch theorem
for curves over finite fields. This Riemann-Roch theorem,
however, can be rephrased in terms of adeles and proved by
means of Fourier analysis. On the other hand, as a
consequence of the explicit formulae, the Riemann
hypothesis can be reformulated as the positivity of the
Weil distribution on the ideles. Thus it should be
possible to give a proof of the Riemann hypothesis via the
Weil positivity that is based entirely on Fourier analysis
on the adeles and ideles. In this spirit, A. Connes
\cite{con} gave an asymptotic trace identity that is
equivalent to the Riemann hypothesis. He managed to give
direct proves for analogous formulae in the local and the
semi-local case, but the global formula can as yet only be
proved as a consequence of the Riemann hypothesis. At this
point it looks as if the difficulty lies in the change of
topology when transgressing from finitely many to
infinitely many places. In this note we restate the trace
formula  as a semi-local trace identity, in which the
global situation makes no appearance at all. The semi-local
trace identity thus becomes equivalent to the Riemann
hypothesis for function fields.

\section{Connes' theorem}
In this section we fix notations and recall Connes' result.
Let $k$ be a global field of positive characteristic $p$.
Then $k$ is the function field of some curve defined over a
finite field. Let $q=p^m$ be the number of elements of the
field of constants in $k$. Let $V$ be the set of
valuations or places of $k$. For each $v\in V$ let $k_v$
be the completion of $k$ at $v$ and let $\CO_v$ be the
ring of integers of the local field $k_v$, i.e. $\CO_v$
consists of all $x\in k_v$ which satisfy $v(x)\ge 0$. For
each $v\in V$ fix a {\it uniformizer at} $p$, i.e. an
element $\pi_v$ of $\CO_v$ such that $v(\pi_v)=1$. Let
$\A$ be the adele ring of $k$, i.e. the subset of the
infinite product $\prod_{v\in V}k_v$ consisting of all
elements $(x_v)_v$ with $x_v\in \CO_v$ for all but
finitely many $v$. We say that $\A$ is the restricted
product of the $k_v$ and write this as $ \A\=
\hat{\prod_v} k_v$. For any subset $S$ of $V$ let
$\A_S=\hat{\prod}_{v\in S}k_v$ and
$\A^S=\hat{\prod}_{v\notin S}k_v$, then
$\A=\A_S\times\A^S$. The ring $\A$ is a locally compact
ring and $k$ embeds diagonally as a discrete subring that
is cocompact as additive group.

For $x\in k_v$ let $|x|_v$ be its modulus, i.e. the unique
positive real number such that for any measurable subset
$A$ of $k_v$ we have $ \mu(xA)\= |x|_v\mu(A)$, where $\mu$
is any additive Haar measure on $k_v$. It then turns out
that $|x|_v=q_v^{-v(x)}$, where $q_v$ is the number of
elements of the  residue class field of $k_v$.

The group of ideles, i.e. the multiplicative group $\A^\times$ of
invertible elements of $\A$ is the restricted product of the
$k_v^\times$ with respect to their compact subgroups
$\CO_v^\times$. In this way the ideles form a locally compact
group whose topology differs from that inherited from the adeles.
The group $k^\times$ embeds diagonally as a discrete subgroup of
$\A^\times$. Let the absolute value on $\A^\times$ be defined as
$|x|=\prod_v |x_v|_v$, which does make sense  since almost all
factors are one. Then this coincides with the modulus of $x$ for
any given additive Haar measure on $\A$. Let $\A^1$ be the set of
all $a\in \A^\times$ such that $|a|=1$, then $k^\times$ forms a
cocompact subgroup of $\A^1$.

The image of $|.|: \A^\times\ra\R$ equals $q_0^\Z$, where
$q_0$ is some positive power of $q$. Let $\pi_V$ be an
element of $\A^\times$ with $|\pi_V| = q_0^{-1}$.

For any subset $S$ of $V$ let $\CO_S=\prod_{v\in S}\CO_v$.
For $S=V$ the ring $\CO_V$ is a compact subring of $\A$ and
$\CO_V^\times$ is a compact subgroup of $\A^\times$. There
is a finite set $\CE\subset \A^1$ such that $\A^1=k^\times
\CE \CO_V^\times$ and the intersection
$k^\times\cap\CO_V^\times$ is the group of nonzero
constants in $k$, i.e. $\F_q^\times$.

Let $\CS(\A)$ be the {\it Schwartz-Bruhat space} of $\A$,
i.e. the space of all locally constant functions on $\A$
with compact support. Any $f\in\CS(\A)$ is a finite sum of
functions of the form $f=\prod_vf_v$, where $f_v$ is
locally constant and of compact support on $k_v$ and
$f(x)=\prod_vf_v(x_v)$; further for all but finitely many
$v$ the function $f_v$ then coincides with $\1_{\CO_v}$,
the characteristic function of $\CO_v\subset k_v$. Fix a
nontrivial additive character $\psi$ on $\A$ which is
trivial on $k$. Note that then the lattice $k$ in $\A$
becomes self-dual \cite{weil}, i.e., for $x\in\A$ we have
$$
\psi(x\ga)\= 1\ \ \ {\rm for\ every\ }\ga\in k\
\Leftrightarrow\ x\in k.
$$
Every additive character on $\A$ decomposes into a product
$$
\psi\=\prod_v\psi_v,
$$
where $\psi_v$ is a character of $k_v$. For $v\in V$ let
$n(v)$ denote the order of $\psi_v$, i.e., $n(v)$ is the
greatest integer $k$ such that $\psi_v$ is trivial on
$\pi_v^{-k}\CO_v$. Then $n(v)=0$ for all but finitely many
$v$ (Cor. 1.IV in \cite{weil}). Let $F_v$ denote the
$\CO_v$-module $\pi_v^{-n(v)}\CO_v$. Then $F_v=\CO_v$ for
all but finitely many $v$.

For any set $S$ of places let $\psi_S=\prod_{v\in
S}\psi_v$. Then $\psi_S$ is a nontrivial additive
character on $\A_S$. Fix a Haar measure $dx$ on $\A_S$ by
the condition that it is self-dual with respect to
$\psi_S$. To explain this, let the Fourier transform on
$\CS(\A_S)$ be defined by
$$
 \hat{f}(x)\= \int_\A f(y)\psi(xy)dy.
$$
Then the measure is normalized so that $\hat{\hat
f}(x)=f(-x)$. More explicitly this means that
$dx=\prod_{v\in S} dx_v$ with $dx_v$ giving the set
$\CO_v$ the volume $q_v^{-\frac{n(v)}2}$.

For any set $S$ of places with at least two elements
define the ring
$$
k_S\= \{ x\in k^\times : |x|_v\le 1\ {\rm for}\ v\notin S\}
$$
and let $k_S^\times$ be the units of this ring, then
$$
k_S^\times\= \{ x\in k^\times : |x|_v=1\ {\rm for}\
v\notin S\}.
$$
Then $k_S^\times$ is a discrete subgroup of $\A_S^\times$,
contained in $\A_S^1$ and such that the quotient
$k_S^\times\bs\A_S^1$ is compact.

For the multiplicative Haar measures we adopt Weil's
normalization as follows: Suppose that $G$ is a locally
compact group and $g\mapsto |g|$ a nontrivial proper
continuous homomorphism to $\R_+^\times$. Then there is a
unique Haar measure on $G$ such that
$$
\vol(\{ g\in G\mid 1\le |g|\le \La\})\ \sim\ \log\La,
$$
as $\La\ra\infty$. We will fix this Haar measure on the
groups $k_v^\times$ for any place $v$ and for
$k_S^\times\bs\A_S^\times$, where $S$ is any set of places
with at least two elements. We will denote these Haar
measures by $d^\times x$. In this normalization we get $
\vol(\CO_v^\times)\= \log q_v,\ \ {\rm and}\ \
\vol(k_S^\times\bs\A_S^1)\=\log q_S, $ where $q_S>1$ is the
generator of the group of absolute values $|x|$, where
$x\in\A_S^\times$. On $\A_S^\times$ we install the measure
given by
$$
\int_{\A_S^\times}\ph(x)\, d^\times
x\=\int_{k_S^\times\bs\A_S^\times}\sum_{\ga\in k_S^\times}
\ph(\ga x)\, d^\times x.
$$

 Note that if $f\in\CS(\A)$ is a product, say,
$f=\prod_vf_v$, then $\hat f=\prod_v\widehat {f_v}$, where
on the right hand side one takes the local Fourier
transforms. Further for a place $v$ and $z\in k_v^\times$
let $\1_{z\CO_v}$ be the characteristic function of the
set $z\CO_v$. Then the local transform satisfies $
\widehat{\1_{z\CO_v}}\=
|z|q_v^{-\frac{n(v)}2}\1_{z^{-1}F_v}. $

Let $ \CS_0(\A)\=\{ f\in\CS(\A) \mid f(0)=0=\hat f(0)\}. $
For $f\in \CS_0(\A)$ let
$$
E(f)(x)\= |x|^{1/2}\sum_{\ga\in k^\times}f(x\ga).
$$
Then $E(f)$ is a function on $\A^\times$, invariant under
$k^\times$, and A. Connes proved in \cite{con} that $E(f)$
is an element of $L^2(k^\times\bs\A^\times)$. Moreover, as
a consequence of Theorem 1 of \cite{con} it follows that
the image of $E$ forms a dense subspace of
$L^2(k^\times\bs\A^\times)$.

Let $\La >0$ be in the value group $q_0^\Z$, and let
$\tilde{Q}_{\La,0}$ be the subspace of $\CS_0(\A)$
consisting of all functions $f$ with $f(x)=0=\hat{f}(x)$
for all $x$ with $|x|>\La$. Let $Q_{\La,0}$ be the closure
in $L^2(k^\times\bs\A^\times)$ of $E(\tilde{Q}_{\La,0})$.
We denote the orthogonal projection onto the space
$Q_{\La,0}$ by the same symbol.

Let $h$ be a Schwartz-Bruhat function on
$k^\times\bs\A^\times$, i.e. the function $h$ is locally
constant and of compact support. We write
$h\in\CS(k^\times\bs\A^\times)$. Given $h$, define an
operator $U(h)$ on $L^2(k^\times\bs\A^\times)$ by
$$
U(h)\ph(x)\=\int_{k^\times\bs\A^\times}|a|^{-\rez{2}}h(a)\ph(xa)d^\times
a.
$$
Let $\La >1$ be in the value group $q_0^\Z$, and let
$$
\log'\La \= \vol(\{a\in k^\times\bs \A^\times | \La^{-1}\le
|a|\le\La\}).
$$
Further let for any $v\in V$: $ h\mapsto
\int'\frac{h(u)}{|u-1|}d^\times u $ be the distribution on
$k_v^\times$ that agrees with $\frac{d^\times u}{|u-1|}$
for $u\ne 1$ and whose Fourier transform vanishes at $1$.
Another way to characterize this distribution is to say it
is the unique distribution that agrees with
$\frac{d^\times u}{|u-1|}$ for $u\ne 1$ and that sends
$h=\1_{\CO_v^\times}$ to zero.

Let $\hat h:\widehat{k^\times\bs\A^\times}\ra\C$ be the
Fourier transform of $h$. In particular write
$$
\hat h(0)\=\int_{k^\times\bs\A^\times}h(x)\ d^\times x,
$$
and
$$
\hat h(1)\=\int_{k^\times\bs\A^\times}h(x)|x|^{-1}\
d^\times x,
$$

The next theorem has been proved by A. Connes in
\cite{con}.

\begin{theorem}
The following are equivalent:
\begin{itemize}
\item[(a)]
For any $h\in\CS(k^\times\bs\A^\times)$ we have, as
$\La\ra\infty$, the asymptotic expansion:
$$
\tr(Q_{\La,0} U(h))\= 2h(1)\log'\La -\hat h(0)-\hat h(1)+
\sum_{v\notin V} \int_{k_v}'\frac{h(u)}{|u-1|}d^\times u
+o(1).
$$
\item[(b)]
All $L$-functions with Gr\"o\ss encharacters on $k$ satisfy the
Riemann hypothesis. \end{itemize}
\end{theorem}

In the light of the fact that $(b)$ is known to be true one
has to read this theorem as follows: $(a)$ is implied by
$(b)$ and an independent proof of $(a)$ gives an
independent proof of $(b)$.

Note that Connes gives a slightly different version of the
trace formula. It is, however, easy to see that the two
versions are equivalent.

\section{The semi-local trace identity}\label{sec 2}
Let $S\subset V$ be a finite set of places that is
supposed to be {\it large enough}. This means that $S$
satisfies the following conditions:
\begin{itemize}
\item $|S|\ge 2$.
\item
The image of the absolute value $|\cdot |_S :
\A_S^\times\ra \R^\times$ is the full value group $q_0^\Z$.
\item
$S$ contains all places $v$ for that the order $n(v)$ of
the character $\psi$ in nonzero.
\item
$S$ is so large, that $c_S=\prod_{v\in S} q_v^{1+n(v)}\ge
1$.
\item
Fix a finite set $\CE\subset\A^1$ such that $\A^1=\CE
k^\times \CO_V^\times$. Then for each $e\in\CE$, the set
$S$ contains all places $v$ with $e_v\notin \CO_v^\times$.
\end{itemize}

If we embed $\A_S^\times$ into $\A^\times$ by $x\mapsto
(x,1,\dots)$ and $S$ is large enough, then the set $\CE$
can be chosen to be contained in $\A_S^\times$. We will
tacitly assume this. We will further assume that $\CE$ is
a set of representatives for $\A^\times /k^\times
O^\times$. Then it also is a set of representatives for
$\A_S^\times /k_S^\times O_S^\times$.

Recall that $k_S^\times$ is a discrete subgroup of
$\A_S^\times$. For $f\in\CS(\A_S)$ let
$$
E_S(f)(x)\= |x|_S^\rez{2} \sum_{\ga\in k_S^\times}f(\ga x).
$$
Then $E_S(f)$ lies in $L^2(k_S^\times\bs\A_S^\times)$, as
is shown in \cite{con}. Further, let
$k_S^0=k_S\setminus\{0\}$ and define
$$
\bar E_S(f)(x)\= |x|_S^\rez{2} \sum_{\ga\in k_S^0}f(\ga x).
$$
We define $\psi_S =\prod_{v\in S}\psi_v$ as a character of
$\A_S$ and the Fourier transform on $\CS(\A_S)$ by $
\hat{f}(x)\=\int_{\A_S}f(y) \psi(xy) dy. $ Let
$\CS_0(\A_S)$ be the space of all $f\in\CS(\A)$ with
$f(0)=0=\hat f(0)$.

\begin{proposition} \label{2.1}
For any $f_S\in\CS_0(\A_S)$ the function $\bar E_S(f_S)$
lies in the space $L^2(k_S^\times\bs \A_S^\times)$.
Furthermore, we have $\bar E_S(f)(x)=\bar E_S(\hat
f)(x^{-1})$.
\end{proposition}

\prf Let $T=V\setminus S$. We extend a given $f_S\in
\CS(\A_S)$ to $f\in \CS(\A)$ by $f=f_S\otimes\1_{\CO_T}$.
A. Connes has shown \cite{con} that $E(f)\in
L^2(k^\times\bs \A^\times)$.

\begin{lemma}
The map
\begin{eqnarray*}
\alpha: k_S^\times\bs\A_S^\times &\ra& k^\times\bs\A^\times
/\CO_T^\times\\
\hspace{10pt} k_S^\times a &\mapsto& k^\times a\CO_T^\times
\end{eqnarray*}
is a homeomorphism, equivariant under the action of
$\A_S^\times$ and Haar-measure preserving, where on
$k^\times\bs\A^\times /\CO_T^\times$ we install the
quotient measure given by the normalized measure (i.e.
$\vol(\CO_T^\times)=1$) on $\CO_T^\times$.
\end{lemma}

\prf The map $\alpha$ is well defined since $k_S^\times
a=k_S^\times a'$ implies $k^\times a=k^\times a'$. It is
clearly continuous and $\A_S^\times$-equivariant. It is
injective since for $a, a'\in \A_S^\times$ embedded into
$\A^\times$ by filling up with ones, the equation
$k^\times a\CO_T^\times = k^\times a'\CO_T^\times$ implies
that there are $x\in k^\times$ and $y\in\CO_T^\times$ such
that $a=xa'y$. Considering the places in $v\in T$ gives
$1=x_vy_v$, which implies $x_v\in\CO_v^\times$ and thus
$x\in k_S^\times$. The map $\alpha$ finally is surjective
since $S$ is large enough and $\A^\times =k^\times \CE
\CO_V^\times \pi_{S}^\Z$, where $\pi_S\in\A_S^\times$ is a
uniformizer, i.e., $|\pi_S|=q_0^{-1}$. The map $\alpha$
preserves Haar measures as a consequence of our
normalizations. Finally, since $\CO_T^\times$ is open in
$\A_T^\times$ it follows that the map $\alpha$ is open,
i.e. the inverse map is continuous. \qed

As a consequence of the lemma, $\alpha$ induces a pullback
isomorphism
$$
\alpha^* : L^2( k^\times \bs \A^\times / \CO_T^\times )
\ra L^2(k_S^\times \bs\A_S^\times).
$$
Now let $x\in\A_S$ and write $(x,1)$ for the element of
$\A$ that coincides with $x$ in $S$ and is equal to $1$
everywhere else. We get $f((x,1))=f_S(x)$, and $
Ef(x,1)\=\bar E_S(f)(x). $ The function $f$ is by
construction invariant under the multiplication by
$\CO_T^\times$, hence lies in the space of invariants $
L^2(k^\times\bs \A^\times)^{\CO_T^\times}\cong
L^2(k^\times \bs \A^\times /\CO_T^\times), $ and we
conclude that $\bar E_S(f_S)=\alpha^*E(F)$. The second
assertion of the proposition follows from Lemma 2 in
Appendix I of \cite{con}.
\qed

Let $\La >0$ and $\tilde{Q}_{S,\La,0}$ be the subspace of
$\CS_0(\A_S)$ consisting of all $f\in \CS_0(\A_S)$ such
that $f(x)=0=\hat{f}(x)$ whenever $|x|>\La$. Let
$Q_{S,\La,0}$ be the closure in
$L^2(k_S^\times\bs\A_S^\times)$ of the space
$E(\tilde{Q}_{S,\La,0})$. Likewise, let $\bar
Q_{S,\Lambda,0}$ be the closure in
$L^2(k_S^\times\bs\A_S^\times)$ of the space $\bar
E_S(\tilde Q_{S,\Lambda,0})$. We also write $
Q_{S,\Lambda,0}$ and $\bar Q_{S,\Lambda,0}$ for the
orthogonal projections.

Let $h\in\CS(k_S^\times\bs\A_S^\times)$, i.e. $h$ is
locally constant and of compact support. Define the
operator $U(h)$ on $L^2(k_S^\times\bs\A_S^\times)$ by
$$
U(h)\ph(x)\=\int_{k_S^\times\bs\A_S^\times}|y|^{-\rez{2}}
h(y)\ph(xy)d^\times y.
$$
Let $\hat h:\widehat{k_S^\times\bs\A_S^\times}\ra\C$ be
the Fourier transform of $h$. In particular write $ \hat
h(0)\=\int_{k_S^\times\bs\A_S^\times}h(x)\ d^\times x, $
and $ \hat
h(1)\=\int_{k_S^\times\bs\A_S^\times}h(x)|x|_S^{-1}\
d^\times x. $ Let $\tilde Q_{S,\La}$ be the subspace of
$\CS(\A_S)$ consisting of all $f\in\CS(\A_S)$ with
$f(x)=0=\hat f(x)$ whenever $|x|>\La$. Let $Q_{S,\La}$ be
the closure in $L^2(k_S^\times\bs\A_S^\times)$ of the
space $E_S(\tilde Q_{S,\La})$.

\begin{theorem}\label{main theorem} (Semi-local trace identity) The following
assertions are equivalent.
\begin{itemize}
\item[(a)] If $h$ is supported in $\{ q_0^{-r}\le |x|\le
q_0^r\}$ and $S$ contains all places $v$ with $q_v\le
q_0^r$, then, as $\La\ra\infty$ we have
$$
\tr Q_{S,\La,0}U(h)\= \tr \bar Q_{S,\La,0}U(h) + o(1).
$$
\item[(b)] Connes' global trace formula.
\item[(c)] The Riemann hypothesis for all $L$-functions
with Gr\"o\ss encharacters.
\end{itemize}
\end{theorem}

Since Connes showed that $(b)$ is equivalent to $(c)$ it
suffices to show that $(a)$ is equivalent to $(b)$. The
proof of this theorem will be given in section \ref{proof
of main}.

\section{A variant of the semi-local trace formula}
Let $S$ be a finite set of places that is large enough
(see section \ref{sec 2}). Let $h\in \CS(k_S^\times\bs
\A_S^\times)$.

\begin{theorem}\label{3.1}
As $\la\ra\infty$, we have
$$
\tr Q_{S,\La,0}U(h)\= 2h(1)\log_S'\La -\hat h(0)-\hat h(1)
+\sum_{v\in S}\int_{k_v^\times}' \frac{h(u)}{|u-1|}
d^\times u +o(1).
$$
\end{theorem}

\prf In \cite{con} it is shown that
$$
\tr(Q_{S,\La}U(h))\= 2h(1)\log'_S\La + \sum_{v\in S}
\int_{k_v^\times}\frac{h(u)}{|u-1|} d^\times u +o(1),
$$
as $\La\ra\infty$. Actually, Connes shows a slightly
different assertion in \cite{con}, namely, instead of one
projection $Q_{S,\La}$ there is a product of two
projections $\hat P_\La P_\La$. But similar to the global
case it is easy to see that, in the absence of infinite
places, Connes' statement is equivalent to the above.

 Let $z=\prod_{v\in S}z_v$ be in $\A_S^\times$. Then,
locally
$$
\widehat{\1_{z_v\CO_v^\times}}\=
|z_v|\,q_v^{-\frac{n(v)}2}\left(\1_{z_v^{-1}F_v}-\rez{q_v}\1_{\pi_v^{-1}z_v^{-1}F_v}\right).
$$
So we get
$$
\widehat{\1_{z\CO_S^\times}}\= |z|\prod_{v\in
S}q_v^{-\frac{n(v)}2}
\left(\1_{z_v^{-1}F_v}-\rez{q_v}\1_{\pi_v^{-1}z_v^{-1}F_v}\right).
$$
Let $c_S=\prod_{v\in S} q_v^{1+n(v)}$. Since $S$ is large
enough we infer that $c_S\ge 1$. Then for $a\in\A_S^\times$
with $\frac{c_S}\La \le |a|\le \La$, the function
$f=\1_{a\CO_S^\times}$ lies in $\tilde Q_{S,\La}$ and $
E_S(f)(x)\= (q-1)|x|^{\rez{2}}\1_{a\CO_S^\times
k_S^\times}(x). $ For $\alpha<\beta$ in the value group
$q_0^\Z$, let $A(\alpha,\beta)$ be a set of
representatives of $ \{ x\in\A_S \mid\alpha \le
|x|\le\beta\}{\big /}\CO_S^\times k_S^\times. $ Note that
$A(\alpha,\beta)$ is a finite set. Let $
f_{1,\La}\=\rez{q-1}\sum_{a\in
A\left(\frac{c_S}\La,\La\right)}\1_{a\CO_S^\times}. $ Then
$f_{1,\La}\in \tilde Q_{S,\La}$ and $ E_S(f_{1,\La})(x)\=
|x|^{\rez 2}\1_{\left\{\frac{c_S}\La \le
|x|\le\La\right\}}. $ Let $f_{0,\La}$ be the Fourier
transform of $f_{1,\La}$.

\begin{lemma}\label{2.5} Suppose that $h$ is supported in
$\{ |x|\le 1\}$. Then
$$
Q_{S,\La}U(h)E_S(f_{1,\La})\ \equiv\
\frac{\frac{c_S}{|b|}-\La^2 q_0}{c_S-\La^2 q_0}\hat h(0)
E_S(f_{1,\La})\ \ \ \mod \left(Q_{S,\La,0}\right).
$$
\end{lemma}

\prf Note first that the assertion only depends on the
$\CO_S^\times$-invariant projection $ h^{\CO_S^\times}(x)\=
\rez{\vol(\CO_S^\times)}\int_{\CO_S^\times}h(xy) d^\times y
$ of $h$. So we may assume that $h$ is
$\CO_S^\times$-invariant. Then it follows that $h$ is a
finite linear combination of functions of the form
$\1_{b\CO_S^\times k_S^\times}$, for $b\in\A_S^\times$. We
may thus assume that $h=(q-1)\1_{b\CO_S^\times
k_S^\times}$ with $|b|\le 1$. Then $ U(h)E_S(f_{1,\La})\=
E_S(R(g)f_{1,\La}), $ where $g=\1_{b\CO_S^\times}$, and $
R(g)f_{1,\La}(x) = \int_{b\CO_S^\times} f_{1,\La}(x y)
d^\times y = {\vol(\CO_S^\times)}f_{1,\La}(xb). $ So that
\eqn
U(h)E_S(f_{1,\La})
&=& \vol(\CO_S^\times)E_S\left( \rez{q-1}\sum_{a\in A\left(\frac{c_S}{\La
|b|},\frac{\La}{|b|}\right)} \1_{a\CO_S^\times}\right).
\endeqn
Since $|b|\le 1$, we get that $Q_{S,\La}U(h)E_S(f_{1,\La})$
equals
$$
\vol(\CO_S^\times)E_S\left(\rez{q-1}\sum_{a\in
A\left(\frac{c_S}{\La |b|},{\La}\right)}
\1_{a\CO_S^\times}\right).
$$
For $f\in \CS(\A_S)$ let $\hat l(f)=\hat
f(0)=\int_{\A_S}f(x) dx$. Then
$$
\hat l\left(\1_{a\CO_S^\times}\right)\=
\int_{\A_S}\1_{a\CO_S^\times}(x) dx\=
|a|\int_{\CO_S^\times} d x.
$$
It follows that there is a $\la\in\C$ with $Q_{S,\La}U(h)
E_S(f_{1,\La})-\la E_D(f_{1,\La})\in Q_{S,\La,0}$ and this
$\la$ is
$$
\la\=\vol(\CO_S^\times)\frac{\sum_{a\in
A\left(\frac{c_S}{\La |b|},\La\right)}|a|}
             {\sum_{a\in A\left(\frac{c_S}{\La},\La\right)}|a|}.
$$
For $\alpha,\beta$ in the value group we can choose $
A(\alpha,\beta)\= \bigcup_{\alpha \le q_0^k\le \beta}
\pi_{S}^{-k} \CE, $ where $\pi_S\in\A_S^\times$ is an
element with $|\pi_S|=q_0^{-1}$. Let $\La=q_0^k$,
$c_S=q_0^{k_0}$, and $v(b)\ge 0$ the valuation of $b$. Then
\eqn
\la &=& \vol(\CO_S^\times)
\frac{|\CE|\sum_{j=k_0-k+v(b)}^k q_0^j}
                        {|\CE|\sum_{j=k_0-k}^k q_0^j}\\
&=& \vol(\CO_S^\times)\frac{\sum_{j=0}^{2k-k_0-v(b)} q_0^j q_0^{k_0-k+v(b)}}
                        {\sum_{j=0}^{2k-k_0} q_0^j q_0^{k_0-k}}\\
&=& \vol(\CO_S^\times)q_0^{v(b)}\frac{1-q_0^{2k-k_0-v(b)+1}}{1-q_0^{2k-k_0+1}}
\= \vol(\CO_S^\times)\rez{|b|}\frac{c_S-\La^2 q_0
|b|}{c_S-\La^2 q_0}.
\endeqn
Since $ \vol(\CO_S^\times k_S^\times /k_S^\times)\=
\frac{\vol(\CO_S^\times)}{q-1}, $ we get $\hat h(0)=
\vol(\CO_S^\times)$. The lemma follows.
\qed

Now write $ f_{1,\La}\= e_{1,\La}+q_{1,\La}, $ where
$q_{1,\La}\in Q_{S,\La,0}$ and $e_{1,\La}$ is orthogonal
to $Q_{S,\La,0}$. Likewise write $ f_{0,\La}\=
e_{0,\La}+q_{0,\La}. $ Let $l(E_S(f))=f(0)$ for $E_S(f)\in
Q_{S,\La}^C$. Note that $l$ is well defined and that $
\ker(l)\= Q_{S,\La,0}^C\oplus \C e_{1,\La}. $ Likewise let
$\hat l(E_S(f))=\hat f(0)$ and note that $ \ker(\hat l)\=
Q_{S,\La,0}\oplus \C e_{0,\La}. $ The operator
$Q_{S,\La}U(h)$ preserves $\ker(l)$ if $\supp h\subset \{
|x|\le 1\}$, and it preserves $\ker(\hat l)$ if $\supp
h\subset\{ |x|\ge 1\}$. Let $ \tilde f_{1\La}\=
\frac{f_{1,\La}}{\norm{E_S(f_{1,\La})}}. $

\begin{lemma} If $h$ is supported in $\{ |x|\ge 1\}$, then,
as $\La\ra\infty$,
$$
Q_{S,\La} U(h) E_S(\tilde f_{1,\La}) \= \hat h(1)
E_S(\tilde f_{1,\La}) +\ph_\La + o(1),
$$
where $\ph_\La\in Q_{S,\La,0}$.
\end{lemma}

\prf Without loss of generality we can assume that
$h=(q-1)\1_{b\CO_S^\times k_S^\times}$ for some $|b|\ge 1$.
Let
$$
f_{1,\La}^1\=\rez{q-1}\sum_{a\in A\left(
\frac{c_S|b|}{\La},\La\right)} \1_{a\CO_S^\times},
$$
and let
$$
\tilde f_{1,\La}^1 \=
\frac{f_{1,\La}^1}{\norm{E_S(f_{1,\La})}},\ \ \tilde
e_{1,\La} \= \frac{e_{1,\La}}{\norm{E_S(f_{1,\La})}},\ \
\tilde q_{1,\La} \=
\frac{q_{1,\La}}{\norm{E_S(f_{1,\La})}}.
$$
Then  $\norm{E_S(\tilde f_{1,\La}-\tilde f_{1,\La}^1)}$
tends to zero as $\La\ra\infty$. Similar to the last proof
we get
$$
U(h) E_S(f_{1,\La}^1)\=\vol(\CO_S^\times) E_S\left(
\rez{q-1}\sum_{a\in A\left(\frac{c_S}\La ,\frac
\La{|b|}\right)} \1_{a\CO_S^\times}\right),
$$
and this also equals $Q_{S,\La}U(h)E_S(f_{1,\La}^1)$.
Repeating the argument of the last lemma we get $
Q_{S,\La} U(h) E_S(f_{1,\La}^1)\ \equiv\ \la
E_S(f_{1,\La}^1), $ with
\eqn
\la &=& \vol(\CO_S^\times)\frac{\sum_{a\in
A(\frac{c_S}{\La},\frac{\La}{|b|})} |a|}
 {\sum_{a\in A(\frac{c_S|b|}{\La},{\La})} |a|}\= \vol(\CO_S^\times)\frac{\sum_{a\in
A(\frac{c_S}{\La},\frac{\La}{|b|})} |a|}
 {\sum_{a\in A\left(\frac{c_S}{\La},\frac{\La}{|b|}\right)} |a b|}\= \rez{|b|}\vol(\CO_S^\times).
\endeqn
Since $\hat h(1)=\rez{|b|}\vol(\CO_S^\times)$, the claim
follows.
\qed

Suppose further that $h$ is supported in $\{ |x|\ge 1\}$.
Then $Q_{S,\La}U(h)$ preserves $\ker(\hat
l)=Q_{S,\La,0}\oplus\C \tilde e_{1,\La}$, and we get
\eqn
Q_{S,\La} U(h) \tilde e_{1,\La} &\equiv & Q_{S,\La} U(h)
E_S(\tilde f_{1,\La})\ \ \ \mod\ \ker(\hat l)\\
&\equiv& \hat h(1) E_S(\tilde f_{1,\La})+ o(1)\ \ \ \mod\
\ker(\hat l)\\
&\equiv& \hat h(1) \tilde e_{1,\La} +o(1)\ \ \ \mod\
\ker(\hat l).
\endeqn
Next we show that $\tilde q_{1,\La}$ tends to zero as
$\la\ra\infty$. We use the fact that for any $\ph\in
Q_{S,\La,0}$ the integral $\int \ph(x) |x|^{\rez 2}
d^\times x$ vanishes, to  compute
\eqn
\sp{\tilde q_{1,\La}}{\tilde q_{1,\La}} &=& \sp {\tilde
q_{1,\La}}{E_S(\tilde f_{1,\La})}\= \int_{|x|\ge
\frac{c_S}\La} \tilde
q_{1,\La}(x)|x|^{\rez 2} d^\times x\\
&=& -\int_{|x| < \frac{c_S}\La} \tilde q_{1,\La}(x)
|x|^{\rez 2} d^\times x.
\endeqn
We infer that
$$
\norm{\tilde q_{1,\La}}^2\ \le\
\sqrt{\int_{|x|<\frac{c_S}{\La}}|\tilde q_{1,\La}(x)|^2
d^\times x\ \int_{|x|<\frac{c_S}{\La}}|x|d^\times x}.
$$
The right hand side tends to zero as $\La\ra\infty$.
Therefore $\tilde q_{1,\La}\ra 0$. Using Lemma \ref{2.5}
we get for $\supp h\subset\{|x|\le 1\}$,
\eqn
Q_{S,\La}U(h) \tilde e_{1,\La} &=& Q_{S,\La} U(h)
E_S(\tilde f_{1,\La}) + o(1)\\
&\equiv& \hat h(0) E_S(\tilde f_{1,\La}) +o(1) \ \ \ \mod
Q_{S,\La,0}\\
&=& \hat h(0) \tilde e_{1,\La} +o(1).
\endeqn
The Fourier transform turns $\ker(\hat
l)=Q_{S,\La,0}\oplus \C \tilde e_{0,\La}$ into
$\ker(l)=Q_{S,\La,0}\oplus \C\tilde e_{1,\La}$. For $\ph\in
L^2(k_S^\times\bs\A_S^\times)$ we have
$\widehat{U(h)\ph}=U(\breve h)\hat\ph$, where $\breve
h(x)=|x|h(x^{-1})$. Since $\hat{\breve h}(0)=\hat h(1)$
and $\hat{\breve h}(1)=\hat h(0)$, we get for
$\supp(h)\subset\{ |x|\ge 1\}$,
$$
Q_{S,\La}U(h)\tilde e_{0,\La}\ \equiv\ \hat h(1)\tilde
e_{0,\La}+o(1) \ \ \ \mod\ker(l).
$$
Likewise, we get for $\supp (h)\subset\{ |x|\ge 1\}$,
$$
Q_{S,\La}U(h)\tilde e_{0,\La}\ \equiv\ \hat h(0)\tilde
e_{0,\La}+o(1) \ \ \ \mod\ker(l).
$$
We have proved the following lemma, which also implies the
theorem.

\begin{lemma}
As $\La\ra\infty$, we have
$$
\tr(Q_{S,\La}-Q_{S,\La,0})U(h)\=\hat h(0) +\hat h(1) +
o(1).
$$
\end{lemma}
\qed

\section{The localization}
Let $h\in\CS(k^\times\bs\A^\times)$. There is a function
$g\in \CS(\A^\times)$ such that $ h(x)=\sum_{\ga\in
k^\times} g(\ga x). $ There is a finite set of places $S$
that is large enough (sec. \ref{sec 2}) such that
$g=g_Sg^S$, where $g_S\in\CS(\A_D^\times)$ and
$g^S=\prod_{v\notin S}\1_{\CO_v^\times}$. For
$x\in\A_S^\times$ set
$$
h_S(x)\=\sum_{\ga\in k_S^\times} g_S(\ga x).
$$

\begin{theorem} \label{4.1}(Localization)
For $\La >1$ we have
$$
\tr Q_{\La,0} U(h)\=\tr \bar Q_{S,\La,0} U(h_S).
$$
\end{theorem}

\prf For $\ph\in L^2(k^\times\bs\A^\times)$ we get $
U(h)\ph\= R(g)\ph, $ where
$$
R(g)\ph(x)\=\int_{\A^\times}g(y)|y|^{-\rez 2} \ph(xy)
d^\times y.
$$
Then $R(g)=R(g_S)R(g^S)$, where, for $\ph\in
L^2(k^\times\bs\A^\times)$ we have
$$
R(g_S)\ph(x)\=
\int_{\A_S^\times}|a|^{-\rez{2}}g_S(a)\ph(xa)d^\times a
$$
and
$$
R(g^S)\ph(x)\=
\int_{\A^{S,\times}}|a|^{-\rez{2}}g^S(a)\ph(xa)d^\times a,
$$
where the measure in the latter integral is the quotient
of the normalized measures on $\A^\times$ and
$\A_S^\times$. Let $T=V\setminus S$. With this measure the
set $\CO_T^\times\subset \A^{S,\times}$ has volume $1$.

It follows that
$$
R(g^S) \= \rez{\vol(\CO_T^\times)}
\int_{\A_T^\times}|a|^{-\rez{2}}\1_{\CO_T^\times}(a)\ph(xa)
d^\times a\=
\rez{\vol(\CO_T^\times)}\int_{\CO_T^\times}\ph(xa) d^\times
a.
$$
Thus it emerges that $R(g^S)$ coincides with $Pr_T$, the
orthogonal projection onto the space $
L^2(k^\times\bs\A^\times)^{\CO_T^\times} $ of
$\CO_T^\times$-invariants. On $L^2(k^\times\bs \A^\times)$
we now have two orthogonal projections, $Q_\La$ and
$R(g^S)$, which commute with each other. So we get
\eqn
\tr Q_{\La,0} U(h) &=& \tr Q_{\La,0} R(g)\= \tr
\left(Q_{\La,0}R(g)\left|\, L^2(k^\times \bs
\A^\times)^{\CO_T^\times}\right.\right)\\
&=& \tr \left( \alpha^* Q_{\La,0}R(g_S)(\alpha^*)^{-1}
\left|\, L^2(k_S^\times\bs\A_S^\times)\right.\right).
\endeqn
It is easy to see that $R(g_S)$ commutes with $\alpha^*$.

\begin{lemma} We have the identity of projections $\alpha^* Q_{\La,0}(\alpha^*)^{-1}=\bar Q_{S\La,0}$.
\end{lemma}

\prf Let $Q_{\La,0}^{\CO_T^\times}$ be the space of
$\CO_T^\times$-invariants in $Q_{\La,0}$. The lemma will
follow from the identity of vector spaces
$\alpha^*\left(Q_{\La,0}^{\CO_T^\times}\right)=\bar
Q_{S,\La,0}$. For this let $f\in\tilde Q_{\La,0}$ and
suppose that $E(f)$ is $\CO_T^\times$-invariant. We may
then assume that $f$ itself is $\CO_T^\times$-invariant.
Then $f$ can be written as a finite sum $f=\sum_j f_j$,
where each $f_j$ lies in $\tilde Q_\La$ and is a product
$f_j=\prod_v f_{j,v}$. We may assume that if $v\notin S$,
then $f_{j,v}=\1_{b\CO_v}$ for some $b\in k_v^\times$. All
but finitely many of the $b$'s can be chosen to be $1$.
Since $S$ is large enough, there is, for each $j$, a
$\ga_j\in k^\times$ such that with
$f_j^{\ga_j}(x)=f_j(\ga_j x)$ we have $
f_{j,v}^{\ga_j}\=\1_{\CO_v} $ for every $v\notin S$. Let $
f_1\=\sum_j f_j^{\ga_j}. $ We claim that $f_1$ lies in
$\tilde Q_{\La,0}$ again. By $f_j\in \tilde Q_\La$ and
$|\ga_j|=1$ we get that $f_j^{\ga_j}\in \tilde Q_\La$. So
we only have to show that $f_1(0)=0=\hat{f_1}(0)$.

For the first, recall that
$$
f_1(0)=\sum_j f_j^{\ga_j}(0)=\sum_j f_j(0)=f(0)=0.
$$
For the second recall
$$
\hat{f_1}(0)\= \sum_j \widehat{f_j^{\ga_j}}(0)\= \sum_j
\hat f_j^{\ga_j^{-1}} (0)\= \sum_f \hat f_j(0)\= \hat
f(0)\= 0.
$$
Finally,
$$
E(f_1)\=\sum_j E(f_j^{\ga_j})\= \sum_jE(f_j)\= E(f),
$$
and this shows that $ E(\tilde Q_{\La,0})^{\CO_T^\times}\=
E(\tilde Q_{\La,0}^T), $ where $\tilde Q_{\La,0}^T$ is the
space of all $f\in\tilde Q_{\La,0}$ that can be written as
a product $f=f_S\left(\prod_{v\notin S}\1_{\CO_v}\right)$
for some $f_S\in \tilde Q_{S,\La,0}$. For such a function
$f$ and $x\in \A_S^\times$ we get
$$
E(f)(x) \= \sum_{\ga\in k^\times} f(\ga x)\= \sum_{\ga\in
k_S^0} f_S(\ga x)\= \bar E_S(f_S)(x).
$$
Since all $f_S$ in $\tilde Q_{S,\La,0}$ can occur, the
lemma follows and this implies the theorem.
\qed

\section{Proof of the main theorem}\label{proof of main}
We will now prove Theorem \ref{main theorem}. We use the
notation of the previous section.

By Theorem \ref{4.1} Connes' trace formula is equivalent to
$$
\tr \bar Q_{S,\La,0} U(h)\= 2h(1)\log'\La -\hat h(0)-\hat
h(1)+ \sum_{v\notin V}
\int_{k_v}'\frac{h(u)}{|u-1|}d^\times u +o(1).
$$
By Theorem \ref{3.1}, on the other hand, we know that
$$
\tr(Q_{S,\La,0}U(h_S))\= 2h_S(1)\log'_S\La -\hat
h_S(0)-\hat h_S(1) + \sum_{v\in S}
\int_{k_v^\times}\frac{h_S(u)}{|u-1|} d^\times u +o(1).
$$
So the main theorem will follow from identifying the right
hand sides of both of theses formulae.

\begin{lemma}
We have
\begin{itemize}
\item[(i)] $h(x)=h_S(x)$ if $x\in\A_S^\times$,
\item[(ii)] $\hat h(0) =\hat h_S(0)$ and $\hat h(1) =\hat
h_S(1)$,
\item[(iii)] $\log_S'\La=\log'\La$.
\end{itemize}
\end{lemma}

\prf For $x\in \A_S^\times$ we have
$$
h(x) \= \sum_{\ga\in k^\times} g(\ga x)\= \sum_{\ga\in
k_S^\times} g_S(\ga x)\= h_S(x).
$$
This proves $(i)$. We further see that
$$
\hat h(0) \= \int_{k^\times\bs\A^\times}h(x)d^\times x\=
\int_{\A^\times}g(x) d^\times x\= \int_{\A_S^\times}
g_S(x) d^\times x\= \hat h_S(0),
$$
as well as $\hat h(1)=\hat h_S(1)$. Finally note that
\eqn
\log_S'\La &=& \int_{k_S^\times\bs\A_S^\times}\sum_{a\in
A(\rez\La,\La)}\1_{a\CO_S^\times k_S^\times}(x) d^\times
x\= \rez{q-1}\sum_{a\in A(\rez\La,\La)} \int_{\A_S^\times}
\1_{a\CO_S^\times}(x) d^\times x\\
&=& \frac{\vol(\CO_S^\times}{q-1}\ \# A\left(\rez \La,\La\right)\= \log'\La.
\endeqn
The last equation follows from the same computation with
$S$ replaced by $V$ and the fact that
$A\left(\rez\La,\La\right)$ also is a set of
representatives of the set of all $x\in \A^\times $ of
absolute value between $\rez\La$ and $\La$ modulo
$k^\times \CO_V^\times$.
\qed

\begin{lemma}
Suppose $h$ is supported in $\{ q^{-r}\le |x|\le q^r\}$.
Then for every $v\in S$ we have
$$
\int_{k_v^\times}'\frac{h_S(u)}{|u-1|} d^\times u\=
\int_{k_v^\times}'\frac{h(u)}{|u-1|} d^\times u.
$$
For $v\notin S$ we have
$\int_{k_v^\times}'\frac{h(u)}{|u-1|} d^\times u=0$.
\end{lemma}

\prf The first assertion follows from $(i)$ of the last
lemma. For the second let $v\notin S$. Then
$$
 \int_{k_v^\times} \frac{h(u)}{|u-1|}
d^\times u\= \sum_{\ga\in k^\times} \int_{k_v^\times}
\frac{g(\ga u)}{|u-1|} d^\times u.
$$
Let $\ga\in k^\times$. If $g(\ga u)\ne 0$ for some $u\in
k_v^\times$, then it follows that $q_0^{-r}\le |\ga|_S\le
q_0^r$ and $|\ga|_w=1$ for every $w\notin S$, $w\ne v$.
Therefore it follows $q_0^{-r}\le |\ga|_v\le q_0^r$ which
implies $|\ga|_v=1$ since $q_v>q_0^r$. Therefore $u\mapsto
g(\ga u)$ is a multiple of $\1_{\CO_v^\times}$ and so
$$
\int_{k_v^\times} \frac{g(\ga u)}{|u-1|} d^\times u \= 0.
$$
\qed

\section{Closing remarks}
In this section we give a reformulation of the semi-local
trace identity.
 Let $C$ be  a compact open subgroup  of
$\A_S^\times$, such that $h$ is invariant under
translations by elements of $C$. Let $P_C$ be the
orthogonal projections onto the space of $C$-invariants.
Then $ P_C\=\rez{\vol(c)}\int_{C} d^\times x. $ This
implies that $P_C$ leaves stable the spaces $Q_{S,\La}$,
$Q_{S,\La,0}$ and $\bar Q_{S,\La,0}$. Further we have that
$U(h)=P_C U(h) P_C$. Let $Q_{S,\La}^C$, $\bar
Q_{S,\La,0}^C$ and $Q_{S,\La,0}^C$ be the subspaces of
$C$-invariants. Every $\ph\in \bar Q_{S,\La,o}^C$ is
supported in $\{\rez\La \le |x|\le\La\}$. Since the set $
\{\rez\La \le |x|\le\La\}/k_S^\times C $ is finite, it
follows that $\bar Q_{S,\La,0}^C$ is finite dimensional. In
particular the set $\bar Q_{S,\La,0}^C$ coincides with
$\bar E_S(\tilde Q_{S,\La,0}^C)$. Let $f\in\tilde
Q_{S,\La,0}$. Then $ \bar E_S(f)(x)\= |x|^{\rez
2}\sum_{\ga\in k_S^0} f(\ga x). $ Let
$R=k_S^0/k_S^\times$. We use the same letter to indicate a
set of representatives for the quotient $R$. Summing over
$k_S^\times$ first and then over $R$ gives the sum
expansion
$$
\bar E_S(f)(x)\=\sum_{r\in R}|r|^{-\rez 2} E_S(f)(rx).
$$
This sum converges pointwise, indeed the sum is locally
finite, but
 it does not converge absolutely in $L^2(k_S^\times \bs
\A_S^\times)$. In this way we get a canonical surjective
map $ T\colon E_S\left(\tilde Q_{S,\La,0}^C\right)\ \ra\
\bar Q_{S,\La,0}^C $ defined by
$$
T(E_S(f))(x) =
\sum_{r\in R}|r|^{-\rez 2} E_S(f)(rx)= \bar E_S(f)(x).
$$
It
is not hard to see that there are coefficients $c_r\in\R$
such that for every $f\in \tilde Q_{S,\La,0}^C$ we have
$$
E_S(f)(x)\= T'(\bar E_S(f))(x)\= \sum_{r\in R} c_r \bar
E_S(f)(rx),
$$
where, again, the sum is locally finite. This implies that
$E(\tilde Q_{S,\La,0}^C)$ is finite dimensional and so
coincides with $Q_{S,\La,0}^C$. Further $T$ is a linear
bijection with inverse $T^{-1}=T'$.

If we extend $T$ to a bicontinuous linear bijection of
$L^2(k_s^\times\bs \A_S^\times)$ that maps the orthogonal
space of $Q_{S,\La,0}$ to the orthogonal space of $\bar
Q_{S,\La,0}$, then we get
$$
T Q_{S,\la,0} T^{-1}=\bar Q_{S,\la,0}.
$$
So that the semi-local trace identity then becomes
$$
\tr Q_{S,\La,0} U(h)\= \tr Q_{S,\La,0} T U(h)T^{-1}.
$$
Note that if $h$ is supported in the norm one elements,
then $U(h)$ leaves invariant the space $Q_{S,\La,0}$ and
it commutes with $T$. So in this case the formula follows
directly. It would be nice to find a direct prove of this
identity in general.

University of Exeter, Mathematics, Exeter EX4 4QE, England


\newpage
\begin{thebibliography}{XXX}

\bibitem{bomb}
 \bf Bombieri, E.:
 \it Counting points on curves over finite fields (d'après S. A. Stepanov).
 \rm Séminaire Bourbaki, 25ème
année (1972/1973), Exp. No. 430, pp. 234--241. Lecture
Notes in Math., Vol. 383, Springer, Berlin, 1974.


\bibitem{con}
 \bf Connes, A.:
 \it Trace formulas in noncommutative geometry and
the zeros of the Riemann zeta function.
 \rm Sel. Math., New Ser. 5, No.1, 29-106 (1999).

\bibitem{con2}
 \bf Connes, A.:
 \it Trace formula on the adèle class space and
Weil positivity.
 \rm Current developments in mathematics, 1997
(Cambridge, MA), 5-64, Int. Press, Boston, MA, 1999.

 \bibitem{con3}
 \bf Connes, A.:
 \it Noncommutative geometry and the Riemann zeta function.
 \rm Mathematics: frontiers and
perspectives, 35-54, Amer. Math. Soc., Providence, RI,
2000.

\bibitem{polya}
 \bf Deitmar, A.:
 \it A Polya-Hilbert operator for automorphic L-function
 \rm to appear in: Indagationes Mathematicae.

\bibitem{weil}
\bf Weil, A: \it Basic number theory. \rm Classics in
Mathematics. Springer-Verlag, Berlin, 1995.

\bibitem{weil2}
 \bf Weil, A.:
 \it Numbers of solutions of equations in finite fields.
 \rm Collected Papers I, 399-410.

\bibitem{weil3}
 \bf Weil, A.:
 \it Sur les Courbes Alg\'ebriques et les Vari\'et\'es s'en
 d\'eduisent.
 \rm Hermann \& Cie, Paris 1948.


\end{thebibliography}
\end{document}